\documentclass{article}
\usepackage{latexsym}
\usepackage{amssymb,amsthm}

\usepackage{amsmath,amsthm,latexsym,amssymb,tikz}
\usetikzlibrary{arrows}

\begin{document}

\newtheorem{defi}{Definition}
\newtheorem{theo}{Theorem}
\newtheorem{coro}{Corollary}
\newtheorem{prop}{Proposition}
\newtheorem{rem}{Rem:}
\newtheorem{la}{Lemma}
\newtheorem{exa}{Example}
\newtheorem{con}{Conjecture}
\newtheorem{que}{Question}

\begin{center}
{\Large Proof of a Conjecture on the Wiener Index of Eulerian Graphs }\\[3mm]
{\large Peter Dankelmann\footnote{Financial support
by the South African National Research Foundation, grant 118521, is
gratefully acknowledged.} (University of Johannesburg)}
\end{center}

\begin{abstract}
The Wiener index of a connected graph is the sum of the distances between
all unordered pairs of vertices. A connected graph is Eulerian if its vertex
degrees are all even. In  
[Gutman, Cruz, Rada, Wiener index of Eulerian Graphs, Discrete Applied Mathematics
132 (2014), 247-250]
the authors proved that the cycle is the unique graph maximising the Wiener index 
among all Eulerian graphs of given order. They also conjectured that for Eulerian 
graphs of order $n \geq 26$ the graph consisting of a cycle on $n-2$ vertices and 
a triangle that share a vertex is the unique Eulerian
graph with second largest Wiener index. The conjecture is known to hold for 
all $n\leq 25$ with exception of six values. 
In this paper we prove the conjecture. 
\end{abstract}
MSC-class: 05C12 (primary) 92E10 (secondary)

\section{Introduction} 

Let $G=(V,E)$ be a finite, connected graph. The Wiener index of $G$ is defined by
\[ W(G) = \sum_{\{u,v\} \subseteq V} d_G(u,v), \]
where $d_G(u,v)$ denotes the usual distance between vertices $u$ and $v$ of $G$, i.e.,
the minimum number of edges on a $(u, v)$- path in $G$. 

The Wiener index, originally conceived by the chemist Wiener \cite{Wie1947},  
has been investigated extensively in the mathematical and chemical literature, 
often under different names, such as transmission, distance, total distance or gross status. 
Several of these results were originally obtained for the closely related 
average distance, also called mean distance, which is defined as 
$\binom{n}{2}^{-1} W(G)$, where $n$ is the order of the graph $G$. 
For some of its chemical applications see, for example, \cite{Rou2002}.  
  
One of the most basic results on the Wiener index states that 
\[ W(G) \leq \binom{n+1}{3} \] 
for every connected
graph on $n$ vertices, and equality holds if and only if $G$ is a path. 

A path has only vertices of degree one and two, so it is reasonable to expect that 
better bounds can be obtained if restrictions are placed on the values of the degrees. 
Upper bounds on the Wiener 
index that take into account not only the order, but also the minimum degree were
given, for example, in \cite{ BeeRieSmi2001, DanEnt2000, KouWin1997}, and it was shown
in \cite{AloDan-manu} that stronger bounds hold in the presence of a vertex of large degree. 
The Wiener index in relation to the inverse degree, i.e., the sum of the inverses of all 
vertex degrees, was considered by Erd\"{o}s, Pach, and Spencer \cite{ErdPacSpe1988}.

Bounds on the Wiener index of trees in terms of vertex degree have also been 
considered extensively. Every tree has minimum degree $1$, so 
it is natural to ask how large or small
the Wiener index can be in trees of given maximum degree. Answering this question for 
the maximum value of the Wiener index is fairly straightforward (see \cite{Ple1984}
and \cite{Ste2008}), however the determination of the minimum Wiener index 
by Fischermann, Hoffmann, Rautenbach, Sz\'{e}kely and Volkmann \cite{FisHofRauSzeVol2002}
required much more effort. For the more general problem of determining the 
extremal values of the Wiener index of a tree with given degree sequence see, for example,  
\cite{CelSchWim2Woe2011}, \cite{SchWagWan2012} and \cite{Wan2008}. 
A good survey of results on the Wiener index of trees before 2000 was given
in \cite{DobEntgut2001}. 

Not the actual value, but the parity of the degrees has been used to bound the Wiener index. 
Trees in which all vertices have odd degree were considered by Lin \cite{Lin2013}, who 
determined their smallest and largest possible Wiener index of trees. 
This result was extended in \cite{ForGutLin2013} with the determination of all such
trees of order $n$ with the largest $\lfloor \frac{n}{4} \rfloor+1$ values of the Wiener index,
see also \cite{For2013}. The smallest and largest Wiener index of a tree whose order and 
number of vertices of even degree are given, was determined in \cite{Lin2014}. 

The Wiener index of connected graphs in which all vertices have even degrees, that is, 
Eulerian graphs, was considered by Gutman, Cruz and Rada \cite{GutCruRad2014}, who 
obtained the following theorem.  

\begin{theo}[Gutman, Cruz and Rada \cite{GutCruRad2014}] 
\label{theo:cycle-has-max-W-among-Eulerian}
Let $G$ be an Eulerian graph of order $n$. Then
\[ W(G) \leq W(C_n), \]
where $C_n$ is the cycle on $n$ vertices. 
Equality holds if and only if $G=C_n$. 
\end{theo}

The authors of Theorem \ref{theo:cycle-has-max-W-among-Eulerian} gave
a direct proof of their result. However, since every Eulerian 
graph is $2$-edge-connected, Theorem 
\ref{theo:cycle-has-max-W-among-Eulerian} can also be obtained as a consequence of
Theoremm \ref{theo:plesnik}(a) below, which states that the cycle is 
the unique graph maximising the Wiener index among all 
$2$-edge-connected graphs of given order.

Gutman, Cruz and Rada \cite{GutCruRad2014} also 
presented a conjecture on the question which Eulerian graph of given
order has the second largest Wiener index. For $n\geq 5$ let $C_{n,3}$ be the graph 
of order $n$ obtained from the disjoint union of two cycles 
on $n-2$ vertices and $3$ vertices, respectively, by identifying 
two vertices, one from each cycle. 
Their conjecture states that $C_{n,3}$ is the unique
graph that has the second largest
Wiener index among all Eulerian graphs of order $n$ for $n\geq 26$. 
It is the aim of this paper to give a proof of this conjecture.  

\begin{figure}
  \begin{center}
\begin{tikzpicture}
  [scale=0.6,inner sep=1mm, 
   vertex/.style={circle,thick,draw}, 
   thickedge/.style={line width=2pt}] 
    \node[vertex] (a1) at (1.5,0) [fill=white] {};
    \node[vertex] (a2) at (4.5,0) [fill=white] {}; 
    \node[vertex] (b1) at (0,1.5) [fill=white] {};
    \node[vertex] (b2) at (3,1.5) [fill=white] {};  
    \node[vertex] (b3) at (6,1.5) [fill=white] {};  
    \node[vertex] (c1) at (1.5,3) [fill=white] {};
    \node[vertex] (c2) at (4.5,3) [fill=white] {};  
     
    \draw[very thick, black] (b2)--(a1)--(b1)--(c1)--(b2)--(a2)--(b3)--(c2)--(b2);  
   
\end{tikzpicture}
\hspace*{1.5em}
\begin{tikzpicture}
  [scale=0.6,inner sep=1mm, 
   vertex/.style={circle,thick,draw}, 
   thickedge/.style={line width=2pt}] 
    \node[vertex] (a1) at (0,0) [fill=white] {};
    \node[vertex] (a2) at (3,0) [fill=white] {}; 
    \node[vertex] (a3) at (6,0) [fill=white] {};     
    \node[vertex] (b1) at (1.5,1.5) [fill=white] {};
    \node[vertex] (b2) at (4.5,1.5) [fill=white] {};  
    \node[vertex] (c1) at (0,3) [fill=white] {};
    \node[vertex] (c2) at (3,3) [fill=white] {};  
    \node[vertex] (c3) at (6,3) [fill=white] {};      
     
    \draw[very thick, black] (b1)--(a1)--(c1)--(b1)--(a2)--(b2)--(c2)--(b1);  
    \draw[very thick, black] (b2)--(a3)--(c3)--(b2);   
\end{tikzpicture}
\hspace*{1.5em}
\begin{tikzpicture}
  [scale=0.6,inner sep=1mm, 
   vertex/.style={circle,thick,draw}, 
   thickedge/.style={line width=2pt}] 
    \node[vertex] (a1) at (1.5,0) [fill=white] {};
    \node[vertex] (a2) at (3,0) [fill=white] {}; 
    \node[vertex] (a3) at (6,0) [fill=white] {};     
    \node[vertex] (b1) at (0,1.5) [fill=white] {};
    \node[vertex] (b2) at (4.5,1.5) [fill=white] {};  
    \node[vertex] (c1) at (1.5,3) [fill=white] {};
    \node[vertex] (c2) at (3,3) [fill=white] {};  
    \node[vertex] (c3) at (6,3) [fill=white] {};      
     
    \draw[very thick, black] (a1)--(a2)--(b2)--(c2)--(c1)--(b1)--(a1);  
    \draw[very thick, black] (b2)--(a3)--(c3)--(b2);   
\end{tikzpicture}
\\[5mm]
\begin{tikzpicture}
  [scale=0.5,inner sep=1mm, 
   vertex/.style={circle,thick,draw}, 
   thickedge/.style={line width=2pt}] 
    \node[vertex] (a1) at (1.5,0) [fill=white] {};
    \node[vertex] (a2) at (4.5,0) [fill=white] {}; 
    \node[vertex] (a3) at (7.5,0) [fill=white] {};     
    \node[vertex] (b1) at (0,1.5) [fill=white] {};
    \node[vertex] (b2) at (3,1.5) [fill=white] {};  
    \node[vertex] (b3) at (6,1.5) [fill=white] {};     
    \node[vertex] (c1) at (1.5,3) [fill=white] {};
    \node[vertex] (c2) at (4.5,3) [fill=white] {};  
    \node[vertex] (c3) at (7.5,3) [fill=white] {};      
     
    \draw[very thick, black] (b1)--(a1)--(b2)--(a2)--(b3)--(a3)--(c3)--(b3)--(c2)--(b2)--(c1)--(b1);    
\end{tikzpicture}
\hspace*{1em}
\begin{tikzpicture}
  [scale=0.5,inner sep=1mm, 
   vertex/.style={circle,thick,draw}, 
   thickedge/.style={line width=2pt}] 
    \node[vertex] (a1) at (1.5,0) [fill=white] {};
    \node[vertex] (a2) at (4.5,0) [fill=white] {}; 
    \node[vertex] (a3) at (7.5,0) [fill=white] {};     
    \node[vertex] (b1) at (0,1.5) [fill=white] {};
    \node[vertex] (b2) at (3,1.5) [fill=white] {};  
    \node[vertex] (b3) at (6,1.5) [fill=white] {};  
    \node[vertex] (b4) at (9,1.5) [fill=white] {};      
    \node[vertex] (c1) at (1.5,3) [fill=white] {};
    \node[vertex] (c2) at (4.5,3) [fill=white] {};  
    \node[vertex] (c3) at (7.5,3) [fill=white] {};      
     
    \draw[very thick, black] (b2)--(a1)--(b1)--(c1)--(b2)--(a2)--(b3)--(c2)--(b2);  
    \draw[very thick, black] (b3)--(a3)--(b4)--(c3)--(b3);   
\end{tikzpicture}
\hspace*{1em}
\begin{tikzpicture}
  [scale=0.5,inner sep=1mm, 
   vertex/.style={circle,thick,draw}, 
   thickedge/.style={line width=2pt}] 
    \node[vertex] (a1) at (0,0) [fill=white] {};
    \node[vertex] (a2) at (3,0) [fill=white] {}; 
    \node[vertex] (a3) at (6,0) [fill=white] {};   
    \node[vertex] (a4) at (9,0) [fill=white] {};         
    \node[vertex] (b1) at (1.5,1.5) [fill=white] {};
    \node[vertex] (b2) at (4.5,1.5) [fill=white] {};  
    \node[vertex] (b3) at (7.5,1.5) [fill=white] {};    
    \node[vertex] (c1) at (0,3) [fill=white] {};
    \node[vertex] (c2) at (3,3) [fill=white] {};  
    \node[vertex] (c3) at (6,3) [fill=white] {};      
    \node[vertex] (c4) at (9,3) [fill=white] {};    
         
    \draw[very thick, black] (a1)--(b1)--(a2)--(b2)--(a3)--(b3)--(a4)--(c4)--(b3)--(c3)--(b2)--(c2)--(b1)--(c1)--(a1);;   
\end{tikzpicture}
\end{center} 
\hspace*{0.1em} \\[-10mm]
  \begin{center}
\begin{tikzpicture}
  [scale=0.5,inner sep=1mm, 
   vertex/.style={circle,thick,draw}, 
   thickedge/.style={line width=2pt}] 
    \node[vertex] (a1) at (2,0) [fill=white] {};
    \node[vertex] (a2) at (4,0) [fill=white] {}; 
    \node[vertex] (b1) at (0.5,0.8) [fill=white] {};
    \node[vertex] (b2) at (5.5,1) [fill=white] {};  
    \node[vertex] (b3) at (7.5,1) [fill=white] {};  
    \node[vertex] (c1up) at (0,3.1) [fill=white] {};
    \node[vertex] (c1down) at (0,1.9) [fill=white] {};
    \node[vertex] (c2) at (6,2.5) [fill=white] {};  
    \node[vertex] (d1) at (0.5,4.2) [fill=white] {};
    \node[vertex] (d2) at (5.5,4) [fill=white] {};
    \node[vertex] (d3) at (7.5,4) [fill=white] {};  
    \node[vertex] (e1) at (2,5) [fill=white] {};
    \node[vertex] (e2) at (4,5) [fill=white] {};
     
    \draw[very thick, black] (a1)--(a2)--(b2)--(c2)--(d2)--(e2)--(e1)--(d1)--(c1up)--(c1down)--(b1)--(a1)
               (c2)--(b3)--(d3)--(c2);                 
\end{tikzpicture}
\hspace*{1em}
\begin{tikzpicture}
  [scale=0.5,inner sep=1mm, 
   vertex/.style={circle,thick,draw}, 
   thickedge/.style={line width=2pt}] 
    \node[vertex] (a1) at (1.5,0) [fill=white] {};
    \node[vertex] (a2) at (4.5,0) [fill=white] {}; 
    \node[vertex] (a3) at (7.5,0) [fill=white] {};   
    \node[vertex] (a4) at (10.5,0) [fill=white] {};         
    \node[vertex] (b1) at (0,1.7) [fill=white] {};
    \node[vertex] (b2) at (3,1.7) [fill=white] {};  
    \node[vertex] (b3) at (6,1.7) [fill=white] {};  
    \node[vertex] (b4) at (9,1.7) [fill=white] {};  
    \node[vertex] (b5) at (12,1.7) [fill=white] {};           
    \node[vertex] (c1) at (1.5,3.4) [fill=white] {};
    \node[vertex] (c2) at (4.5,3.4) [fill=white] {};  
    \node[vertex] (c3) at (7.5,3.4) [fill=white] {};     
    \node[vertex] (c4) at (10.5,3.4) [fill=white] {};          
     
    \draw[very thick, black] (b2)--(a1)--(b1)--(c1)--(b2)--(a2)--(b3)--(c2)--(b2);  
    \draw[very thick, black] (b3)--(a3)--(b4)--(c3)--(b3); 
    \draw[very thick, black] (b4)--(a4)--(b5)--(c4)--(b4);       
\end{tikzpicture}
\caption{The Eulerian graphs of second largest Wiener index for $n=7, 8, 9, 10, 11, 13$.}
\label{fig:exceptions}
\end{center}
\end{figure}
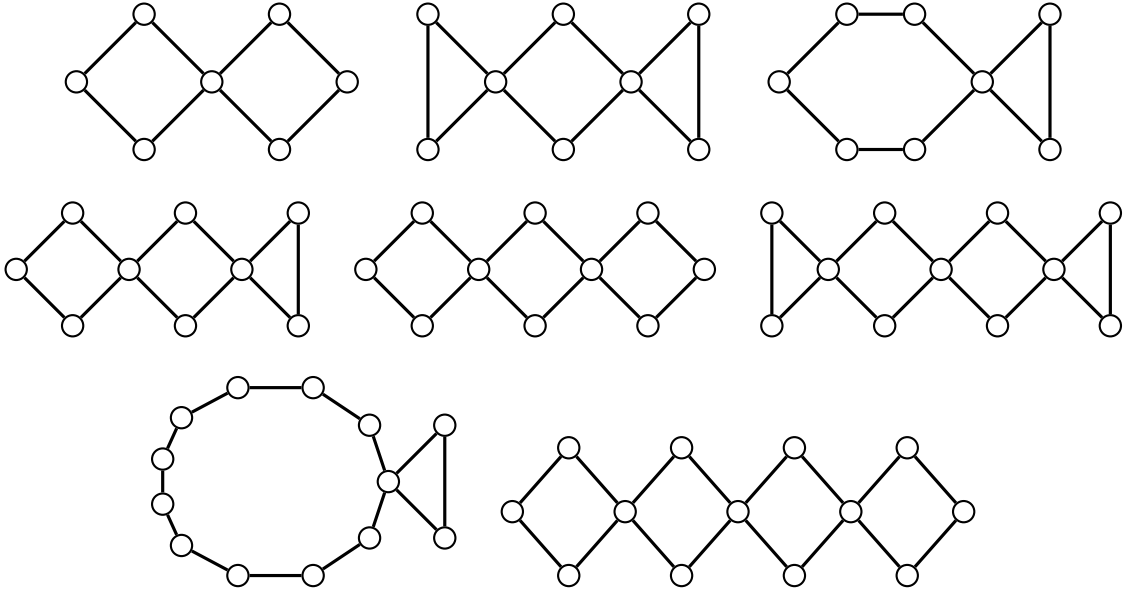

It was verified in  \cite{GutCruRad2014} that the conjecture holds for 
all values of $n$ up to $25$ except $n \in \{7, 9 \}$, 
for which there are other extremal graphs of larger Wiener
index than $C_{n,3}$, and $n \in \{8, 10, 11, 13\}$, for which
$C_{n,3}$ has second largest Wiener index, but there exists another
graph of the same Wiener index. All Eulerian graphs of order 
$7, 8, 9, 10, 11, 13$ that have the second largest Wiener index
are shown in Figure \ref{fig:exceptions}. 

The main result of this paper reads as follows. 

\begin{theo}  \label{theo:Eulerian-with-2nd-largest-W}
Let $G$ be an Eulerian graph of order $n$ with $n\geq 26$ that is not 
a cycle. Then
\[ W(G) \leq W(C_{n,3}) \]
Equality holds if and only if $G=C_{n,3}$. 
\end{theo}

The determination of the unique Eulerian graph with second largest Wiener index 
is  reminiscent of the corresponding problem for $2$-connected graphs. 
A {\em cutvertex} of a connected graph is a vertex whose removal disconnects the graph. 
A connected graph with no cutvertex is said to be $2$-{\em connected}, and a 
{\em block} of a graph is a maximal subgraph that is $2$-connected. 
Plesn\'\i k \cite{Ple1984} proved that among all $2$-connected graphs of given order, the cycle
is the unique graph maximising the Wiener index. While the proof of this result is relatively 
straightforward, determining the $2$-connected graph of given order with the second largest 
Wiener index requires significantly more effort, see the paper by
Bessy, Dross, Knor and \u{S}krekovski \cite{BesDroKnoSkr2020}. The proof of 
Theorem \ref{theo:Eulerian-with-2nd-largest-W} given in the present paper suggest that the 
situation is no different for Eulerian graphs. 
We note that Plesn\'\i k's result on $2$-connected graphs was, asymptotically, extended to  
$k$-connected graphs in \cite{DanMukSwa2009}. 

We note also a certain analogy between our result, and the determination of the largest 
Wiener index of a connected graph with given order and number of cutvertices 
in \cite{BesDroKnoSkr2020-2} and \cite{BesDroHriKnoSkr2020}. If the number of 
cutvertices is sufficiently small relative to the order, then the extremal graph
consists of a path, whose one end is attached to a cycle. So among all graphs
with exactly one cutvertex, the extremal graph has two blocks whose order is 
as unequal as possible. Similarly, the extremal graph $C_{n,3}$ has two blocks,
which are as unequal as possible, given the restriction that every block of an Eulerian graph
has at least three vertices.


The notation we use is as follows. If $G$ is a graph, then $V(G)$ and $E(G)$
denote the vertex set and the edge set, respectively, of $G$. The {\em order} 
$n(G)$ and the {\em size} $m(G)$ are the number of vertices and edges, respectively, of $G$. 
If $G$ and $H$ are graphs with $V(H) \subseteq V(G)$ and $E(H) \subseteq E(G)$,
then we say that $H$ is a {\em subgraph} of $G$ and write $H \leq G$.
If $A \subseteq V(G)$, then $G[A]$ denotes the subgraph of $G$ {\em induced} 
by $A$, i.e., the graph whose vertex set is $A$, and whose edges are exactly
the edges of $G$ joining two vertices of $A$. 
 
If $v$ is a vertex of $G$, then $N_G(v)$ denotes the {\em neighbourhood} of 
$v$, i.e., the set of all vertices of $G$ adjacent to $v$. For $i \in \mathbb{N}$ 
we define the {\em $i$-th neighbourhood} of $v$, $N_i(v)$, to be the set 
of vertices at distance exactly $i$ from $v$, and we let $n_i(v) = |N_i(v)|$.  
The {\em degree} of $v$ in $G$, i.e., the value $n_1(v)$, is denoted by ${\rm deg}_G(v)$. 

A {\em cutset} of $G$ is a set $S \subseteq V$  such that $G-S$, the graph
obtained from deleting all vertices in $S$ and all edges incident with
vertices in $S$ from $G$, is disconnected. 
An {\em edge-cut} of $G$ is a set $E_1 \subseteq E(G)$ such that $G-E_1$,
the graph obtained from $G$ by deleting all edges in $E_1$, is disconnected. 
Let $k\in \mathbb{N}$.  
We say that $G$ is {\em $k$-connected ($k$-edge-connected)} if $G$ contains
no cutset (no edge-cut) with fewer than $k$ elements. 
A {\em cutvertex} is a vertex
$v$ with the property that $\{v\}$ is a cutset. 
An endblock a graph $G$ is a block of $G$ that contains only one cutvertex. 
It is known that every connected graph that is not $2$-connected has at least
two endblocks.

If $S$ is a cutset of $G$ and $H$ a component of $G-S$, then
we say that $G[V(H)\cup S]$ is a {\em branch} of $G$ at $S$. If 
$S=\{v\}$, then we say that $H$ is a branch at $v$.

The {\em total distance} of a vertex $v$, $\sigma_G(v)$, is defined as 
the sum $\sum_{y \in V(G)} d_G(v,y)$. 
By $\sigma_G(A)$ we mean $\sum_{y \in V(G)-A} d_G(y,A)$, where 
$d_G(y,A)$ is defined as $\min_{a\in A} d_G(y,a)$.  

The {\em eccentricity} $e(v)$ of a vertex $v$ of $G$ is the distance from 
$v$ to a vertex farthest from $v$ in $G$.

\section{Preliminary Results}

In this section we present definitions and results that will be used in the proof 
of Theorem \ref{theo:Eulerian-with-2nd-largest-W}. 
We begin with some bounds on the Wiener index and on the total distance 
of vertices in $2$-connected and $2$-edge-connected graphs.

\begin{theo}[Plesn\'\i k \cite{Ple1984}]
  \label{theo:plesnik}
(a) Let $G$ be a $2$-edge-connected graph of order $n$. Then 
\[ W(G) \leq  \left\{ \begin{array}{cc}
      \frac{n^3}{8} & \textrm{if $n$ is even,} \\
      \frac{n^3-n}{8} & \textrm{if $n$ is odd.}
       \end{array} \right.  \] 
Equality holds if and only if $G$ is a cycle. \\       
(b) Let $G$ be a $2$-connected graph of order $n$ and $v$ a vertex of $G$. 
Then
\[ \sigma_G(v) \leq  \left\{ \begin{array}{cc}
      \frac{n^2}{4} & \textrm{if $n$ is even,} \\
      \frac{n^2-1}{4} & \textrm{if $n$ is odd.}
       \end{array} \right.  \] 
Equality holds if $G$ is a cycle. \\       
(c) Let $G$ be a $2$-edge-connected graph of order $n$ 
and $v$ a vertex of $G$. Then
\[  \sigma_G(v) \leq   \frac{n(n-1)}{3}. \]      
\end{theo}

\begin{coro}  \label{la:W-and-sigma-largest-in-cycle}
Let $G$ be a $2$-connected graph of order $n$ and $u,w$ two vertices of $G$.
Let $u_1, u_2$ be two adjacent vertices of the cycle $C_n$. 
Then 
\[ \sigma_G(\{u,w\})  
      \leq \sigma_{C_n}(\{u_1,u_2\}).  \] 
\end{coro}

{\bf Proof:}
Let $G'$ be the $2$-connected graph obtained from $G$ by adding a new
vertex $z$ and joining it to $u$ and $w$. Then 
\[ \sigma(z,G') 
             = \sum_{x\in V(G)} \big(1+d_G(x,\{u,w\})\big) 
             = n + \sigma_{G}(\{u,w\}).   \] 
Let $C_n'$ be the graph obtained from $C_n$ by adding 
adding a new vertex $y$ and joining it to two adjacent vertices $u_1$ and $u_2$
of $C_n$. As above,  
\[ \sigma(y,C_n') 
             = \sum_{x\in V(C_n)} \big(1+d_{C_n}(x,\{u_1,u_2\})\big) 
             = n + \sigma_{C_n}(\{u_1,u_2\}).   \] 
Clearly,
removing the edge $u_1u_2$ from $C_n'$ does not change $\sigma(y)$. 
But $C_n'-u_1u_2$ is  $C_{n+1}$, so by Theorem \ref{theo:plesnik}(b),
we have $\sigma(z,G') \leq \sigma(y,C_n')$, which implies the statement 
of the lemma. \hfill $\Box$

  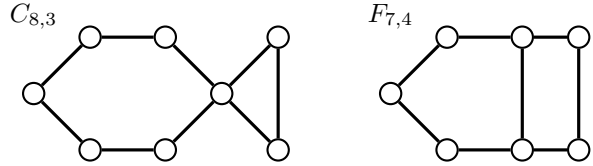
\begin{figure}[h]
  \begin{center}
\begin{tikzpicture}
  [scale=0.5,inner sep=1mm, 
   vertex/.style={circle,thick,draw}, 
   thickedge/.style={line width=2pt}] 
    \begin{scope}[>=triangle 45]    
     \node[vertex] (a1) at (1.5,0) [fill=white] {};
     \node[vertex] (a2) at (3.5,0) [fill=white] {};
     \node[vertex] (a3) at (6.5,0) [fill=white] {};
     \node[vertex] (b1) at (0,1.5) [fill=white] {};
     \node[vertex] (b2) at (5,1.5) [fill=white] {};
     \node[vertex] (c1) at (1.5,3) [fill=white] {};
     \node[vertex] (c2) at (3.5,3) [fill=white] {};
     \node[vertex] (c3) at (6.5,3) [fill=white] {};        
    \draw[very thick] (a1)--(a2)--(b2)--(c2)--(c1)--(b1)--(a1)  
          (b2)--(a3)--(c3)--(b2); 
   \node [above] at (0,3) {$C_{8,3}$};         
    \end{scope}
\end{tikzpicture}   
\hspace*{2em} 
\begin{tikzpicture}
  [scale=0.5,inner sep=1mm, 
   vertex/.style={circle,thick,draw}, 
   thickedge/.style={line width=2pt}] 
    \begin{scope}[>=triangle 45]    
     \node[vertex] (a1) at (1.5,0) [fill=white] {};
     \node[vertex] (a2) at (3.5,0) [fill=white] {};
     \node[vertex] (a3) at (5,0) [fill=white] {};
     \node[vertex] (b1) at (0,1.5) [fill=white] {};
     \node[vertex] (c1) at (1.5,3) [fill=white] {};
     \node[vertex] (c2) at (3.5,3) [fill=white] {};
     \node[vertex] (c3) at (5,3) [fill=white] {};        
    \draw[very thick] (a1)--(a2)--(a3)--(c3)--(c2)--(c1)--(b1)--(a1)  (c2)--(a2); 
   \node [above] at (0,3) {$F_{7,4}$};         
    \end{scope}
\end{tikzpicture} 
\hspace*{2em} 
\caption{Graphs defined in Definitions  \ref{defi:C(n,a)}.}
\label{fig:examples}
\end{center}
\end{figure}

\begin{defi}   \label{defi:C(n,a)}
(a) Let $n, a \in \mathbb{N}$ with $3 \leq a \leq n-2$. 
Then $C_{n,a}$ denotes the graph of order $n$ obtained from two disjoint
cycles $C_{a}$ and $C_{n+1-a}$ by identifying a vertex of $C_{a}$ 
with a vertex of $C_{n+1-a}$. \\
(b) Let $n, a \in \mathbb{N}$ with $3 \leq a \leq n-1$. 
Then $F_{n,a}$ denotes the graph of order $n$ obtained from two disjoint
cycles $C_{a}$ and $C_{n+2-a}$ by choosing two adjacent vertices 
$u, v$ of $C_{a}$ and two adjacent vertices $u', v'$ of $C_{n+2-a}$  
and identifying $u$ with $u'$ and $v$ with $v'$. 
\end{defi}

The Wiener index of the graph $C_{n,a}$ was evaluated in 
\cite{GutCruRad2014}. Specifically for $C_{n,3}$ we have 
\begin{equation} 
W(C_{n,3}) = \left\{ \begin{array}{cc} 
\frac{1}{8}n^3- \frac{1}{4}n^2 + \frac{3}{2}n - 2 & \textrm{if $n$ is even,} \\
\frac{1}{8}n^3- \frac{1}{4}n^2 + \frac{11}{8}n - \frac{9}{4} & \textrm{if $n$ is odd.} 
\end{array} \right. \label{eq:value-W(C(n,3))} 
\end{equation}
In our proofs below we make use of the fact that 
$\frac{1}{8}n^3- \frac{1}{4}n^2 + \frac{11}{8}n - \frac{9}{4} 
\leq W(C_{n,3})
\leq \frac{1}{8}n^3- \frac{1}{4}n^2 + \frac{3}{2}n - 2$ 
for all $n \in \mathbb{N}$ with $n\geq 5$, irrespective of the parity of $n$.

\begin{la} {\rm (Gutman, Cruz, Rada \cite{GutCruRad2014}}  \label{la:Gutman}
If $n \in \mathbb{N}$ is even, $n \geq 6$, then 
\[ W(C_{n,3}) > W(C_{n,4}) > \ldots > W(C_{n,n/2-1})  > W(C_{n,n/2}). \] 
If $n \in \mathbb{N}$ is odd, $n \geq 11$ and $n= 4k+3$ for some $k\in \mathbb{N}$, then 
\[ W(C_{n,3}) > W(C_{n,4}) > \ldots > W(C_{n,2k}) > W(C_{n,2k+2})  > W(C_{n,2k+1}). \] 
If $n \in \mathbb{N}$ is odd, $n \geq 11$ and $n= 4k+1$ for some $k\in \mathbb{N}$, then 
\[ W(C_{n,3}) > W(C_{n,4}) > \ldots > W(C_{n,2k-2}) > W(C_{n,2k}) > W(C_{n,2k-1})  > W(C_{n,2k+1}). \]
For $n=7,9$ we have $W(C_{7,4}) > W(C_{7,3})$ and $W(C_{9,4}) > W(C_{9,3}) > W(C_{9,5})$. 
\end{la}

\begin{la}   \label{la:C(n,3)>F(n,a)}
Let $n\geq 26$ and $4 \leq a \leq n-2$. Then 
\[ W(F_{n,a}) \leq W(C_{n,3}), \]
with equality only if $a=4$ or $a=n-2$. 
\end{la}

{\bf Proof:}
A tedious but straightforward calculation yields that 
\[ W(F_{n,a}) = \left\{ \begin{array}{cc}
\frac{1}{8} \big[ a(n-2)(a-n-2) + n(n^2 + 2n -4) \big] & \textrm{if $n$ even, $a$ even,} \\
\frac{1}{8} \big[ a(n-2)(a-n-2) + n(n^2 + 2n -4) -3n+6 \big] & \textrm{if $n$ even, $a$ odd,} \\
\frac{1}{8} \big[ a(n-2)(a-n-2) + n(n^2 + 2n -4) -n-a+2 \big] & \textrm{if $n$ odd, $a$ even,} \\
\frac{1}{8} \big[ a(n-2)(a-n-2) + n(n^2 + 2n -4) + a -2n \big] & \textrm{if $n$ odd, $a$ odd.} 
\end{array} \right. \]
Since $F_{n,a}=F_{n,n+2-a}$, we may assume that $a \leq \lfloor \frac{n+2}{2} \rfloor$. 
The derivative with respect to $a$ of the four terms on the right hand side above
equals $\frac{1}{8}((n-2)(2a-n-2))$ if $n$ is even, 
$\frac{1}{8}((n-2)(2a-n-2)-1)$ if $n$ is odd and $a$ is even, and
$\frac{1}{8}((n-2)(2a-n-2)+1)$ if $n$ is even and $a$ is odd. 
Hence each of these four terms is strictly decreasing in $a$. It thus follows that 
$W(F_{n,a}) \leq W(F_{n,4})$ if $a$ is even, with equality only if $a=4$, and 
$W(F_{n,a}) \leq W(F_{n,5})$ if $a$ is odd, with equality only if $a=5$. 
Now an easy calculation shows that $W(F_{n,5}) < W(F_{n,4}) = W(C_{n,3})$. 
Hence the lemma follows. \hfill $\Box$

\begin{coro} \label{coro:cycle-plus-triangle}
Let $G$ be a graph of order $n\geq 26$ obtained from a cycle $C_n$ by adding 
three edges between vertices of $C_n$ that are not in $E(C_n)$ and
that form a triangle.  Then $W(G) < W(C_{n,3})$. 
\end{coro}

{\bf Proof:}
Let the three edges added to $C_n$ be $uv, vw, wu$. Then for at least one of 
these three edges, $uv$ say, we have $C_n+uv = F_{n,a}$ for some $a$ with
$4 \leq a \leq \frac{n+1}{2}$. Applying Lemma \ref{la:Gutman} yields that 
$W(G) < W(F_{n,a}) \leq W(C_{n,3})$. \hfill $\Box$

\section{Excluding $2$-Connected Counterexamples}

The goal of this section is to prove that an Eulerian graph of given 
order that is not a cycle, and which has maximum Wiener index among such 
graphs, cannot be $2$-connected. 
Hence it will suffice to prove Theorem \ref{theo:Eulerian-with-2nd-largest-W} for
graphs that have a cutvertex.

We begin by showing that Theorem \ref{theo:Eulerian-with-2nd-largest-W} holds for graphs that are
obtained from two $2$-connected graphs by gluing them together at two
vertices, provided one of them contains a spanning cycle. 

\begin{la}  \label{la:if-branches-are-cycles}
(a) Let $G$ be a $2$-connected graph of order $n$.  If $G$ contains a cutset $\{u,w\}$ 
with 
the property that the union of some, but not all, branches of $G$ at $\{u,w\}$ 
has exactly  $a$ vertices and contains a spanning cycle, and the union of the 
remaining branches is $2$-connected, then
\[ W(G) \leq W(F_{n,a}). \]
Equality implies that $G=F_{n,a}$. \\
(b) If, in addition, $G$ is Eulerian, $n\geq 26$ and $4 \leq a\leq n-2$,  
then $W(G) <  W(C_{n,3})$.  
\end{la}

{\bf Proof:}
Let $A$ be the vertex set of the union of the branches at $\{u,w\}$ that contains 
a spanning cycle, and let $B$ be the vertex set of the union of 
the remaining branches. Let $a=|A|$ and $b=|B|$. Then  $A\cap B =\{u,w\}$. 
Let $H=G[B]$ and let $C$ be a spanning cycle of $G[A]$. 
We denote the set of vertices $x$ of $A-\{u,w\}$ for which $d_C(u,x) < d_C(w,x)$ 
($d_C(u,x) > d_C(w,x)$, $d_C(u,x)=d_C(w,x)$) by $U$ ($W$, $S$).  Then 
\begin{eqnarray}
W(G) & = & W_G(A) + W_G(B) 
      - d_G(u,w) 
     + \sum_{ x \in U \cup W \cup S, y \in B -\{u,w\}} d_G(x,y) 
      \nonumber \\
& \leq & W(C) + W(H) - d_G(u,w) 
    + \sum_{x\in U, y \in V(H)-\{u,w\}} \big(d_C(x,u) + d_H(u,y)\big)   \nonumber \\
   & & +  \sum_{x\in W, y \in V(H)-\{u,w\}} \hspace*{-2em} \big(d_C(x,w) + d_H(w,y)\big) 
     + \sum_{x\in S, y \in V(H)-\{u,w\}} \hspace*{-2em} \big(d_C(x,\{u,w\}) + d_H(\{u,w\},y)\big)  \nonumber \\
& = &  W(C)  + W(H) - d_G(u,w) +  (b-2) \sigma_C(u,U) 
    + |U| (\sigma_H(u)-d_H(u,w))   \nonumber \\ 
   & &  +  (b-2) \sigma_C(w,W)  + |W| (\sigma_H(w)-d_H(u,w))      
    +  (b-2) \sigma_C(\{u,w\},S)   \nonumber  \\
    & &  +|S| \sigma_H(\{u,w\}).     \label{eq:one-branch-a-cycle-1}
\end{eqnarray}
Let $C_a$ and $C_b$ be the cycles of the graph $F_{a,b}$ defined above, and let
$u'$ and $w'$ be the two adjacent vertices of $F_{a,b}$ shared by $C_a$ and $C_b$. 
Let $U'$ ($W'$, $S'$) be the set of 
vertices $x$ of $C_a-\{u,w\}$ with $d(u',x) < d(w',x)$ ($d(u',x) > d(w',x)$, $d(u',x)=d(w',x)$).
As above, we have
\begin{eqnarray}
W(F_{a,b}) 
& = &  W(C_a)  + W(C_b) - d_{F_{a,b}}(u',w') +  (b-2) \sigma_{C_a}(u',U') 
    + |U'| (\sigma_{C_b}(u')-d_{F_{a,b}}(u',w'))   \nonumber \\ 
   & &  +  (b-2) \sigma_{C_a}(w',W')  + |W'| (\sigma_{C_b}(w')-d_{F_{a,b}}(u',w'))      
    +  (b-2) \sigma_{C_a}(\{u',w'\},S') \nonumber \\
    & &  +|S'| \sigma_{C_b}(\{u',w'\}).
      \label{eq:one-branch-a-cycle-2}
\end{eqnarray}
Since $C$ and $C_a$ are cycles, we have
$|U|=|W|$ and $|U'|=|W'|$. 
Subtracting \eqref{eq:one-branch-a-cycle-2} from 
\eqref{eq:one-branch-a-cycle-1} yields thus
\begin{eqnarray*}
W(F_{a,b}) - W(G) & \geq & 
    \big(W(C_a) - W(C)\big) + \big(W(C_b) - W(H)\big) 
    + \big( d_G(u,w) - d_{F_{a,b}}(u',w')\big) \\
  & &  + (b-2) \big[\sigma_{C_a}(u',U') +  \sigma_{C_a}(w',W') +  \sigma_{C_a}(\{u',w'\},S')
             -  \sigma_{C}(u,U) \\
  & &        -  \sigma_{C}(w,W) -  \sigma_{C}(\{u,w\},S)  \big]
    + |U'| \big[ \sigma_{C_b}(u')  + \sigma_{C_b}(w') -  2d_{F_{a,b}}(u',w') \big] \\
  & &  + |S'| \sigma_{C_b}(\{u',w'\})
    - |U| \big[ \sigma_{H}(u)  + \sigma_{H}(w) -  2d_{H}(u,w) \big]
    - |S| \sigma_{H}(\{u,w\})    
\end{eqnarray*}
We now argue that the right hand side of the last inequality is nonnegative. 
Clearly, $C$ and $C_a$ are isomorphic, so $W(C)-W(C_a)=0$.
Since $H$ is $2$-connected, and thus $2$-edge-connected, we have 
$W(C_b)-W(H)\geq 0$ by Theorem \ref{theo:plesnik}(a).  
Since $d_{F_{a,b}}(u',w')=1$ we have $d_G(u,w) - d_{F_{a,b}}(u',w') \geq 0$. 
Also
 $\sigma_{C_a}(u',U') +  \sigma_{C_a}(w',W') +  \sigma_{C_a}(\{u',w'\},S') 
  -  \sigma_{C}(u,U) - \sigma_{C}(w,W) -  \sigma_{C}(\{u,w\},S)
= \sigma_{C_a}(\{u'w'\})  -  \sigma_{C}(\{u,w\})$, but 
$\sigma_{C_a}(\{u'w'\})  -  \sigma_{C}(\{u,w\}) \geq 0$ by 
Lemma \ref{la:W-and-sigma-largest-in-cycle}. \\
We now bound the remaining expression,
$ |U'| \big( \sigma_{C_b}(u')  + \sigma_{C_b}(w') -  2d_{C_b}(u',w') \big) 
  + |S'| \sigma_{C_b}(\{u',w'\})
    - |U| \big( \sigma_{H}(u)  + \sigma_{H}(w) -  2d_{H}(u,w) \big)
    - |S| \sigma_{H}(\{u,w\})$, 
which we denote by $f$. In order to complete the proof of the lemma
it remains to show that $f\geq 0$. 

We have  $a=2|U'|+|S'| = 2|U|+|S|$. 
In $C_a$, vertices $u'$ and $w'$ 
are adjacent, so there is exactly one vertex equidistant from $u'$ and $w'$ if 
$a$ is odd, and there is no vertex equidistant from $u'$ and $w'$ if $a$ is even.
Hence $|S'|=1$ if $a$ is odd, and $|S'|=0$ if $a$ is even. 
In $C$ the vertices $u$ and $w$ are not necessarily adjacent, so we have  
$|S|=1$ if $a$ is odd, and $|S|\in \{0,2\}$ if $a$ is even. 
We conclude that if $a$ is odd, then $|U|=|U'|$ and $|S|=|S'|$, and
if $a$ is even then either $|U|=|U'|$ and $|S|=|S'|$, or 
$|U'|=|U|+1$, $|S'|=0$, and $|S|=2$.
If  $|U'|=|U|$ and $|S|=|S|$, then 
$f = |U|\big(  \sigma_{C_b}(u') - \sigma_{H}(u) + \sigma_{C_b}(w') 
   -  \sigma_{H}(w) +  2d_{H}(u,w) -  2d_{C_b}(u',w') \big) 
   +  |S| \big( \sigma_{C_b}(\{u',w'\} -  \sigma_{H}(\{u,w\}) \big)$. 
Each of the terms $|U|$, $|S|$, $ \sigma_{C_b}(u') - \sigma_{H}(u)$, 
$\sigma_{C_b}(w')-  \sigma_{H}(w)$, $2d_{H}(u,w) -  2d_{C_b}(u',w')$, 
and $\sigma_{C_b}(\{u',w'\} -  \sigma_{H}(\{u,w\})$ is nonnegative,
hence $f\geq 0$ in this case. 
If  $|U'|=|U|+1$ and $|S|=2$, $|S'|=0$, then 
$f =  |U| \big( \sigma_{C_b}(u') -  \sigma_{H}(u)  + \sigma_{C_b}(w') 
         - \sigma_{H}(w) + 2d_{H}(u,w) -  2d_{C_b}(u',w') \big)
        + \sigma_{C_b}(u')  + \sigma_{C_b}(w') -  2d_{C_b}(u',w') 
        - 2\sigma_{H}(\{u,w\})$. 
As above, each of the terms  $|U|$,  $ \sigma_{C_b}(u') - \sigma_{H}(u)$, 
$\sigma_{C_b}(w')-  \sigma_{H}(w)$, $2d_{H}(u,w) -  2d_{C_b}(u',w')$, 
is nonnegative.  We also have 
$\sigma_{C_b}(u')  + \sigma_{C_b}(w') -  2d_{C_b}(u',w')  
        - 2\sigma_{H}(\{u,w\}) \geq 0$
since 
$\sigma_{C_b}(u')  + \sigma_{C_b}(w') -  2d_{C_b}(u',w') = 
\sum_{x\in V(C_b)-\{u',w'\} } \big(d_{C_b}(u',x) + d_{C_b}(w',x)\big)
  \geq  \sum_{x\in V(C_b)-\{u',w'\} } 2 \min \{d_{C_b}(u',x), d_{C_b}(w',x)\}
  = 2\sigma_{C_b}(\{u',w'\})$.
Hence $f\geq 0$ also in this case. This proves the desired 
bound on $W(G)$.  

Now assume that $W(G) = W(F_{n,a})$. Then we have equality between the
corresponding terms in 
\eqref{eq:one-branch-a-cycle-1} and \eqref{eq:one-branch-a-cycle-2},
in particular $W(G[A]) = W(C_a)$ and $W(H)= W(C_b)$.
This implies by Theorem \ref{theo:plesnik}(a) that $G[A]$ and $H$ are cycles
of length $a$ and $b$, respectively. 
We also have $d_C(u,w)=1$. It follows that $G=F_{n,a}$. \\
(b) If $G$ is Eulerian, then $G\neq F_{n,a}$ and so 
$W(G) < W(F_{n,a})$. By Lemma \ref{la:C(n,3)>F(n,a)} we have $W(F_{n,a}) \leq W(C_{n,3})$,
and (b) follows. 
\hfill $\Box$

\begin{la}   \label{la:excluding-2-connected-counterex}
Let $n\in \mathbb{N}$ with $n\geq 26$. Among all Eulerian graphs of order $n$ that are not cycles, 
let $G$ be one that has maximum Wiener index. Then $G$ has a cutvertex. 
\end{la}

{\bf Proof:}
Suppose to the contrary that $G$ is $2$-connected. We first prove that
\begin{equation}  \label{eq:triangle-has-vertex-of-degree-2}
\textrm{every triangle of $G$ contains a vertex of degree $2$.}
\end{equation}
Suppose to the contrary that $G$ contains a triangle $u_1 u_2 u_3$ with 
${\rm deg}(u_i) >2$ for $i=1,2,3$. Let $E'$ be the edge set of 
this triangle. Then $G-E'$ is connected since otherwise, if $G-E'$ is
disconnected, the vertices $u_1, u_2$ and $u_3$ are not all in the same 
component of $G-E'$, so there exists a component of $G-E'$ containing only one
vertex, $u_1$ say, of the triangle. This implies that $u_1$ is a cutvertex
of $G$, a contradiction to $G$ being $2$-connected. Hence 
$G-E'$ is connected. Clearly, $G-E'$ is also Eulerian, and 
$W(G-E') > W(G)$. 
By our choice of $G$, the graph $G-E'$ is a cycle. But then $G$ is obtained from a cycle by
adding the edges of a triangle, and so $W(G) < W(C_{n,3})$ by 
Corollary \ref{coro:cycle-plus-triangle}. 
This contradicts the choice of $G$ as having
maximum Wiener index, and so
\eqref{eq:triangle-has-vertex-of-degree-2}
follows.

  \begin{figure}[h]
  \begin{center}
\begin{tikzpicture}
  [scale=0.5,inner sep=1mm, 
   vertex/.style={circle,thick,draw}, 
   thickedge/.style={line width=2pt}] 
    \begin{scope}[>=triangle 45]    
     \node[vertex] (a1) at (1,0) [fill=white] {};
     \node[vertex] (a2) at (2.5,0) [fill=white] {};
     \node[vertex] (a3) at (4,0) [fill=white] {};
     \node[vertex] (a4) at (5.5,0) [fill=white] {};     
     \node[vertex] (b1) at (0,1) [fill=white] {};
     \node[vertex] (c1) at (1,2) [fill=white] {};
     \node[vertex] (c2) at (2.5,2) [fill=white] {};
     \node[vertex] (c3) at (4,2) [fill=white] {};      
     \node[vertex] (c4) at (5.5,2) [fill=white] {};             
    \draw[very thick] (a1)--(b1)--(c1);
    \draw[very thick] (a2)--(a3)--(a4)  (c2)--(c3)--(c4);    
    \draw[very thick] (c4)--(a3)--(c3)--(a4);     
    \draw[very thick, dotted] (a1)--(a2)  (c1)--(c2);           
   \node [above] at (4,2.2) {$\overline{v}$}; 
   \node [left] at (-0.2,1) {$v$};          
   \node [below] at (4,-0.2) {$w$};     
   \node [above] at (5.5,2.2) {$u_1$};    
   \node [below] at (5.5,-0.2) {$u_2$};              
    \end{scope}
\end{tikzpicture}   
\hspace*{0.3em} 
\begin{tikzpicture}
  [scale=0.5,inner sep=1mm, 
   vertex/.style={circle,thick,draw}, 
   thickedge/.style={line width=2pt}] 
    \begin{scope}[>=triangle 45]    
  
     \node[vertex] (a1) at (1,0) [fill=white] {};
     \node[vertex] (a2) at (2.5,0) [fill=white] {};
     \node[vertex] (b1) at (0,1) [fill=white] {};
     \node[vertex] (b2) at (2.5,1) [fill=white] {};     
     \node[vertex] (c1) at (1,2) [fill=white] {};
     \node[vertex] (c2) at (2.5,2) [fill=white] {};
    \draw[very thick] (a1)--(b1)--(c1)--(a1);
    \draw[very thick] (a1)--(a2) (a1)--(b2)--(c1)  (c2)--(c1);    
    \draw[very thick] (3.0,-0.3)--(a2)--(3.0,0.3) 
    (3.0,1.7)--(c2)--(3.0,2.3);  
   \node [left] at (-0.2,1) {$v$};      
   \node [above] at (1,2.2) {$\overline{v}$};    
   \node [below] at (1,-0.2) {$w$};     
   \node [above] at (2.5,2.2) {$u_1$};   
   \node [right] at (2.7,1.0) {$u_2$};       
   \node [below] at (2.5,-0.2) {$u_3$};             
    \end{scope}
\end{tikzpicture}   
\hspace*{0.3em} 
\begin{tikzpicture}
  [scale=0.5,inner sep=1mm, 
   vertex/.style={circle,thick,draw}, 
   thickedge/.style={line width=2pt}] 
    \begin{scope}[>=triangle 45]    
     \node[vertex] (a1) at (1,0) [fill=white] {};
     \node[vertex] (a2) at (2.5,0) [fill=white] {};
     \node[vertex] (a3) at (4,0) [fill=white] {};
     \node[vertex] (a4) at (5.5,0) [fill=white] {};     
     \node[vertex] (b1) at (0,1) [fill=white] {};
     \node[vertex] (b2) at (5.5,1) [fill=white] {};     
     \node[vertex] (c1) at (1,2) [fill=white] {};
     \node[vertex] (c2) at (2.5,2) [fill=white] {};
     \node[vertex] (c3) at (4,2) [fill=white] {};      
     \node[vertex] (c4) at (5.5,2) [fill=white] {};             
    \draw[very thick] (a1)--(b1)--(c1);
    \draw[very thick] (a2)--(a3)--(a4)  (c2)--(c3)--(c4);    
    \draw[very thick] (c3)--(b2)  (c3)--(a4);     
    \draw[very thick, dotted] (a1)--(a2)  (c1)--(c2);           
    \draw[very thick] (6.2,1.4)--(b2)--(6.2,0.6) (6.2,1.6)--(c4)--(6.2,2.4) 
                (6.2,-0.4)--(a4)--(6.2,0.4);  
   \node [left] at (-0.2,1) {$v$};                  
   \node [above] at (4,2.2) {$\overline{v}$};    
   \node [below] at (4,-0.2) {$w$};     
   \node [above] at (5.5,2.2) {$u_1$};    
   \node [below] at (6.35,1.35) {$u_2$};     
  \node [below] at (5.5,-0.2) {$u_3$};                 
    \end{scope}
\end{tikzpicture}   
\hspace*{0.3em}
\begin{tikzpicture}
  [scale=0.5,inner sep=1mm, 
   vertex/.style={circle,thick,draw}, 
   thickedge/.style={line width=2pt}] 
    \begin{scope}[>=triangle 45]    
     \node[vertex] (a1) at (1,0) [fill=white] {};
     \node[vertex] (a2) at (2.5,0) [fill=white] {};
     \node[vertex] (a3) at (4,0) [fill=white] {};
     \node[vertex] (a4) at (5.5,0) [fill=white] {};     
     \node[vertex] (b1) at (0,1) [fill=white] {};
     \node[vertex] (b2) at (5.5,1) [fill=white] {};     
     \node[vertex] (c1) at (1,2) [fill=white] {};
     \node[vertex] (c2) at (2.5,2) [fill=white] {};
     \node[vertex] (c3) at (4,2) [fill=white] {};      
     \node[vertex] (c4) at (5.5,2) [fill=white] {};             
    \draw[very thick] (a1)--(b1)--(c1);
    \draw[very thick] (a2)--(a3)--(a4)  (c2)--(c3)--(c4);    
    \draw[very thick] (c3)--(b2)  (c3)--(a4)   (c4)--(a3)--(b2);     
    \draw[very thick, dotted] (a1)--(a2)  (c1)--(c2);           
    \draw[very thick] (6.2,1.4)--(b2)--(6.2,0.6) (6.2,1.6)--(c4)--(6.2,2.4) 
                (6.2,-0.4)--(a4)--(6.2,0.4);  
   \node [left] at (-0.2,1) {$v$};                  
   \node [above] at (4,2.2) {$\overline{v}$};    
   \node [below] at (4,-0.2) {$w$};     
   \node [above] at (5.5,2.2) {$u_1$};    
   \node [below] at (6.35,1.35) {$u_2$};     
  \node [below] at (5.5,-0.2) {$u_3$};                 
    \end{scope}
\end{tikzpicture}   
\caption{Cases 1, 2A, 2B, and 2C in the proof of Lemma \ref{la:excluding-2-connected-counterex}.}
\label{fig:cases-in-lemma-4}
\end{center}
\end{figure}
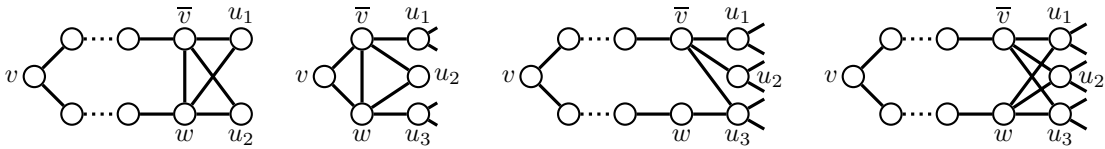

Since $G$ is not a cycle, it has a vertex of degree greater than $2$. 
For $v\in V(G)$ we define $\overline{v}$ to be a nearest vertex of 
degree greater than $2$ (with ties broken arbitrarily) and let $f(v)=d(v,\overline{v})$. 
Note that $v=\overline{v}$ if and only if ${\rm deg}(v) >2$. 
Since $G$ is Eulerian, we have ${\rm deg}(\overline{v}) \geq 4$ for
every $v\in V$. 
For $v\in V$, we have $\overline{v} \in N_{f(v)}(v)$, and thus 
$N(\overline{v}) \subseteq N_{f(v)-1}(v) \cup N_{f(v)}(v) \cup N_{f(v)+1}(v)$.  
We claim that 
\begin{equation}  \label{eq:neighbours-of-v-bar}
\textrm{If ${\rm deg}(v)=2$, then} \ |N(\overline{v}) \cap N_{f(v)-1}(v)| =1.
\end{equation}
Indeed, if ${\rm deg}(v) =2$ then the neighbour, $w$ say, of $\overline{v}$ 
on a shortest $(v,\overline{v})$-path is in $N_{f(v)-1}(v)$.  
If there was a second neighbour $w'$ of $\overline{v}$ in $N_{f(v)-1}(v)$, 
then the vertices in $\bigcup_{i=0}^{f(v)-1} N_i(v) \cup \{\overline{v}\}$ 
would induce a cycle as a subgraph whose only vertex of degree greater than
$2$ in $G$ is $\overline{v}$, implying that $\overline{v}$ is a cutvertex of 
$G$, contradicting the $2$-connectedness of $G$. This proves
\eqref{eq:neighbours-of-v-bar}. \\
Let $v$ be a vertex of degree $2$.  
By the definition of  $f(v)$, all vertices in $\bigcup_{i=0}^{f(v)-1} N_i(v)$ have 
degree $2$. It is easy to see that, since $G$ is $2$-connected, this implies
\begin{equation} \label{eq:first-f(v)-layers-have-2-vertices}
n_1(v) = n_2(v) =\cdots = n_{f(v)}(v) =2. 
\end{equation} 
Let $N_{f(v)}(v)=\{\overline{v}, w\}$. It follows from \eqref{eq:neighbours-of-v-bar} that  
$n_{f(v)}(v) + n_{f(v)+1}(v) \geq {\rm deg}(\overline{v}) \geq 4$, so
$n_{f(v)+1}(v) \geq 2$. We consider three cases, depending on the value 
$n_{f(v)+1}(v)$. \\[1mm]
{\sc Case 1:} There exists $v\in V(G)$ with $n_{f(v)+1}(v) = 2$. \\
Let $N_{f(v)+1}(v)=\{u_1, u_2\}$. Since ${\rm deg}_G(\overline{v}) \geq 4$ it  
follows that $\overline{v}$ is adjacent to $u_1, u_2$, $w$ and a vertex in $N_{f(v)-1}(v)$,
we have ${\rm deg}_G(\overline{v}) = 4$. 
Since $w$ is adjacent to a vertex in $N_{f(v)+1}(v)$, otherwise 
$\overline{v}$ would be a cutvertex, to a vertex in $N_{f(v)-1}(v)$, and
also to $\overline{v}$, it follows that ${\rm deg}_G(w)>2$, and thus 
${\rm deg}_G(w)=4$,  so $w$ is also adjacent to $u_1$ and $u_2$.

Now $\overline{v}$, $w$ and $u_i$ form a triangle for $i=1,2$. Since
${\rm deg}_G(\overline{v}) = {\rm deg}_G(w) =4$, it follows by  
\eqref{eq:triangle-has-vertex-of-degree-2} that 
${\rm deg}_G(u_1) = {\rm deg}_G(u_2) =2$.
So the vertices in $N_{f(v)+1}(v)$ have only neighbours in 
$N_{f(v)+1}(v) \cup N_{f(v)}(v)$. This implies that $e_G(v) = f(v)+1$, and so
$V(G)=\bigcup_{i=0}^{f(v)+1} N_i(v)$. It follows that $G$ consists of the cycle
induced by $\bigcup_{i=0}^{f(v)} N_i(v)$ and the two additional vertices 
$u_1$ and $u_2$ of degree two, both adjacent to $\overline{v}$ and $w$. 
Hence $G$ is the first graph depicted in Figure \ref{fig:cases-in-lemma-4}.
Applying Lemma \ref{la:if-branches-are-cycles} to the cutset $\{\overline{v},w\}$ 
now yields that $W(G)<W(C_{n,3})$. 
This contradiction to the maximality of $W(G)$ proves the lemma in Case 1.  \\[1mm]
{\sc Case 2:} There exists $v\in V(G)$ with $n_{f(v)+1} =3$. \\
Let $N_{f(v)+1}(v)=\{u_1, u_2, u_3\}$. We consider subcases as follows. \\[1mm]
{\sc Case 2a:} $\overline{v}w \in E(G)$. \\
The set $\{\overline{v},w\}$ is a cutset. Its branch containing $v$
is a cycle of length $2f(v)+1$, and the union of the other branches is $2$-connected since
$\overline{v}$ and $w$ are adjacent. Hence we have 
$W(G) \leq W(F_{n,2f(v)+1})$ by Lemma \ref{la:if-branches-are-cycles}(a).
If $f(v)\geq 2$, then $4\leq 2f(v)+1 \leq n-2$ and so 
$W(F_{n,2f(v)+1}) < W(C_{n,3})$ by 
Lemma \ref{la:if-branches-are-cycles}(b), a contradiction to the 
maximality of $W(G)$. Hence $f(v)=1$ and $\overline{v} \in N_1(v)$. 

Then ${\rm deg}(w)>2$, since otherwise $w$ would not have neighbours in $N_{2}(v)$, 
and $\overline{v}$ would be a cutvertex, a contradiction. Hence ${\rm deg}(w)\geq 4$. From 
\eqref{eq:neighbours-of-v-bar} and $n_{1}(v) + n_{2}(v) = 5$ we 
conclude that both, $\overline{v}$ and $w$, have degree $4$, and both are
adjacent to exactly two vertices in $N_{2}(v)$.  
We may assume that $\overline{v}$ is adjacent to $u_1$ and $u_2$, while 
$w$ is adjacent to $u_2$ and $u_3$. 
The situation is depicted in the second graph of Figure \ref{fig:cases-in-lemma-4}. 
Since $\overline{v}$, $w$ and $u_2$ form
a triangle, we have ${\rm deg}(u_2)=2$ by \eqref{eq:triangle-has-vertex-of-degree-2}. 
This implies that $\{\overline{v},w\}$ is a cutset of adjacent
vertices with at least three branches, and the branches containing 
$v$ and $u_2$ are $3$-cycles, so their union contains a spanning cycle, whose
length is $4$. Since
$\overline{v}$ and $w$ are adjacent, the union of the remaining branches
is $2$-connected. Hence it follows by 
Lemma \ref{la:if-branches-are-cycles} that $W(G) < W(C_{n,3})$.
This contradiction to the maximality of $W(G)$ proves the lemma
in this case.  \\[1mm]
{\sc Case 2b:} $\overline{v}w \notin E(G)$ and ${\rm deg}(w)=2$. \\
Then $w$ has a unique neighbour, $u_3$ say, in $N_{f(v)+1}(v)$. 
The set $\bigcup_{i=0}^{f(v)} N_i(v) \cup \{u_3\}$ induces a cycle in
$G$ in which only $\overline{w}$ and possibly $u_3$ have degree
greater than $2$. 
The situation is depicted in the third graph of Figure \ref{fig:cases-in-lemma-4}. 
The set $\{\overline{v}, u_3\}$ is a cutset.
The branch containing $v$ and $w$ induces a cycle of length $2f(v)+2$, 
and the union of the remaining branches is $2$-connected 
since $\overline{v}$ and $u_3$ are adjacent. Since $4\leq 2f(v)+2 \leq n-2$, we 
have $W(G) < W(C_{n,3})$ by Lemma \ref{la:if-branches-are-cycles}, 
a contradiction to the maximality of $W(G)$.  \\[1mm]
{\sc Case 2c:} $\overline{v}w \notin E(G)$ and ${\rm deg}(w)>2$. \\
By \eqref{eq:neighbours-of-v-bar} and $n_{f(v)+1}(v)=3$ if follows
that $\overline{v}$ and $w$ are both adjacent to $u_1$, $u_2$ and $u_3$. 
If at least one vertex in $\{u_1, u_2, u_3\}$, $u_1$ say, has 
degree $2$, then the union of the two branches at the cutset 
$\{ \overline{v}, w\}$ containing $v$ and $u_1$ has at least four vertices and a 
spanning cycle, while the union of the remaining branches is $2$-connected. 
Hence it follows from Lemma \ref{la:if-branches-are-cycles} that 
$W(G) < W(C_{n,3})$, contradicting the maximality of $W(G)$. 
So we may assume that $u_1$, $u_2$ and $u_3$ all have degree greater
than $2$. 
The situation is depicted in the fourth graph of Figure \ref{fig:cases-in-lemma-4}. 
Let $E'$ be the edge set of the $4$-cycle $\overline{v}, u_1, w, u_2, \overline{v}$.
Then $G-E'$ is connected since otherwise, similarly to the proof of 
\eqref{eq:triangle-has-vertex-of-degree-2}, one of the vertices 
$u_1$, $u_2$ or $u_3$ would be a cutvertex of $G$. 
Since all vertices in $G-E'$ have even degree, 
it follows that $G-E'$ is Eulerian. Since at least one vertex of $G-E'$ has
degree greater than $2$, viz $u_3$, we conclude that $G-E'$ is not
a cycle. But $W(G-E') > W(G)$, a contradiction to the
maximality of $W(G)$. \\[1mm]
{\sc Case 3:} $n_{f(v)+1}(v) \geq 4$ for all $v\in V$. \\
Let $v\in V(G)$ be fixed. 
We first show that 
\begin{equation} \label{eq:bound-on-d(v)-and-f(v)}
\sigma_G(v) \leq \left\{ \begin{array}{cc}
\frac{1}{4}n^2 - n + \frac{11}{4} + 2f(v) & \textrm{if $n$ is odd,}\\
\frac{1}{4}n^2 - n + 3 + 2f(v) & \textrm{if $n$ is even}
\end{array} \right. 
\end{equation} 
We note that $n_0(v)=1$, 
$n_1(v)=n_2(v)= \ldots = n_{f(v)}(v)=2$, and $n_{f(v)+1}(v)\geq 4$ imply that
$n\geq 5+2f(v)$, so $f(v) \leq \lfloor \frac{n-5}{2} \rfloor$. 
Let $k=e(v)$. 
Then $\sigma_G(v) = \sum_{i=0}^k in_i(v)$. The values $n_i(v)$ satisfy the following 
conditions: (i) $n_0(v)=1$ and (ii) $\sum_{i=0}^k n_i(v)=n$.  
Since $G$ is $2$-connected, we have (iii) $n_i(v) \geq 2$ for $i=1,2,\ldots,k-1$,
and (iv) $n_{f(v)+1}(v) \geq 4$ by the defining condition
of Case 3. 

In order to bound $\sum_{i=1}^k in_i(v)$ from above, assume that $n$ and $f(v)$ are fixed, and that 
integers $k, n_0, n_1,\ldots,n_k$ are chosen to maximise $\sum_{i=1}^k in_i$ subject to
conditions (i)-(iv). Then  
$n_0=1$, and $n_i =2$ for all $i\in \{1,2,\ldots,k-1\} - \{f(v)+1\}$, 
since otherwise, if $n_i>2$, we can modify the sequence $n_0,\ldots,n_k$ by 
decreasing $n_i$ by $1$ and increasing $n_{i+1}$ by $1$ to obtain a new
sequence $n_0',\ldots,n_k'$ which satisfies (i)-(iv), but for 
which $\sum_{i=0}^k in_i' > \sum_{i=0}^k in_i$, a contradiction. 
The same argument yields that $n_{f(v)+1}=4$, and also that 
$n_k\in \{1,2\}$ if $k>f(v)+1$. Therefore, 
if $n$ is odd we have $k=\frac{n-3}{2}$ and 
$\sum_{i=0}^k in_i = \frac{1}{4}n^2 - n + \frac{11}{4} + 2f(v)$, 
and if $n$ is even we have $k=\frac{n-2}{2}$, $n_k=1$ and 
$\sum_{i=0}^k in_i = \frac{1}{4}n^2 - n + 30 + 2f(v)$,  
which is \eqref{eq:bound-on-d(v)-and-f(v)}. 

Summation of \eqref{eq:bound-on-d(v)-and-f(v)} over all $v\in V(G)$
yields 
\begin{equation} \label{eq:bound-on=W-and-f(v)}
2W(G) = \sum_{v\in V(G)} \sigma_G(v)   
    \leq  \left\{ \begin{array}{cc} 
    \frac{1}{4}n^3 - n^2 + \frac{11}{4}n + 2\sum_{v\in V(G)} f(v) & \textrm{if $n$ is odd,}\\
\frac{1}{4}n^3 - n^2 + 3n + 2\sum_{v\in V(G)} f(v) & \textrm{if $n$ is even.}
\end{array} \right.   
\end{equation}
We now bound $\sum_{v \in V(G)}f(v)$. 
Since $G$ is an Eulerian graph but not a cycle, $G$ contains a vertex 
$w$ of degree at least $4$. Since for $i\in \{0,1,\ldots,e(w)\}$ every vertex
$v \in N_i(w)$ satisfies $f(v) \leq d(v,w)$, we have 
\[ \sum_{v\in V(G)} f(v) \leq \sum_{v\in V(G)} d(v,w) = \sigma_G(w).  \]
Now $G$ has more than one vertex of degree greater than two, since 
otherwise such a vertex would be a cutvertex, contradicting the 
$2$-connectedness of $G$. That implies that the strict inequality
$\sum_{v\in V(G)} f(v) < \sigma_G(w)$ holds. 
Noting that $f(w)=0$, we obtain by \eqref{eq:bound-on-d(v)-and-f(v)} that
\begin{equation} \label{eq:bound-on-sum-of-f(v)}
\sum_{v\in V(G)} f(v) < \sigma_G(w) \leq \left\{ \begin{array}{cc}
\frac{1}{4}n^2 - n + \frac{11}{4}  & \textrm{if $n$ is odd,}\\
\frac{1}{4}n^2 - n + 3  & \textrm{if $n$ is even}.
\end{array} \right. 
\end{equation}
From \eqref{eq:bound-on=W-and-f(v)} and \eqref{eq:bound-on-sum-of-f(v)}
we get 
\[ W(G) < \left\{ \begin{array}{cc} 
    \frac{1}{8}n^3  - \frac{1}{4}n^2 + \frac{3}{8}n + \frac{11}{4} & \textrm{if $n$ is odd,}\\
\frac{1}{8}n^3 - \frac{1}{4}n^2 + \frac{1}{2}n + 3  & \textrm{if $n$ is even.} 
\end{array} 
\right.  \]
But the right hand side of the last inequality equals $W(C_{n,3})$.  This contradiction
to the maximality of $W(G)$ completes the proof.  \hfill $\Box$

\section{Completing the proof of Theorem \ref{theo:Eulerian-with-2nd-largest-W} }

{\bf Proof of Theorem \ref{theo:Eulerian-with-2nd-largest-W}:}
Suppose to the contrary that the theorem is false, and let $n$ be the smallest value
with $n\geq 26$ for which the theorem fails. 
Let $G$ be an Eulerian gaph of order $n$ that is not a cycle, and that has maximum 
Wiener index among all such graphs. By Lemma \ref{la:excluding-2-connected-counterex},
$G$ has a cutvertex, so $G$ is not $2$-connected. 
Then $G$ hat at least two endblocks. Let $H$ be a smallest 
endblock of $G$, let $v$ be the cutvertex 
of $G$ contained in $H$, and let $K$ be the union of the branches at $v$ distinct from $H$. 
Let $A$ and $B$ be the vertex set of $H$ and $K$, respectively, and let  $a=|A|$ and $b=|B|$. 
Then $b=n-a+1$, and since $H$ is a smallest endblock we have $a \leq \frac{n+1}{2}$. 
We have 
\begin{eqnarray}
W(G) & = & \sum_{\{x,y\} \subseteq A} d_{H}(x,y) + \sum_{\{x,y\} \subseteq B} d_{K}(x,y) 
    + \sum_{x\in A-\{v\}} \sum_{y \in B-\{v\}} \big( d_{H}(x,v) + d_{K}(v,y) \big) 
    \nonumber \\
&=&  W(H) + W(K) + (a-1) \sigma_{K}(v) + (b-1) \sigma_{H}(v).  \label{eq:split-W-into-2}
\end{eqnarray}
Since $H$ is an endblock, $H$ is $2$-connected, but $K$ may or may not be $2$-connected. 
 \\[1mm]
{\sc Case 1:} $K$ is $2$-connected. \\
Similarly to \eqref{eq:split-W-into-2} we obtain for the graph $C_{n,a}$ and its two 
blocks $C_a$ and $C_{b}$ that  
\[ W(C_{n,a}) 
  =  W(C_{a}) + W(C_{b}) + (a-1) \sigma_{C_b}(w) + (b-1)\sigma_{C_a}(w), \]
where $w$ is the cutvertex of $C_{n,a}$.
Since $H$ and $K$ are $2$-connected, we have by Theorem \ref{theo:plesnik} 
that  $W(H) \leq W(C_a)$, $W(K) \leq W(C_b)$, 
$\sigma_{K}(v) \leq \sigma_{C_b}(w)$ and $\sigma_{H}(v) \leq \sigma_{C_a}(w)$.  
Hence we have $W(G) \leq W(C_{n,a})$. 
By Lemma \ref{la:Gutman} we have $W(C_{n,a}) \leq W(C_{n,3})$, and so
we have $W(G) \leq W(C_{n,3})$, as desired. \\
Assume that $W(G) = W(C_{n,3})$. Then  
$W(K) = W(C_b)$, and so $K=C_b$, and similarly $H=C_a$
by Theorem \ref{theo:plesnik}. Now 
Lemma \ref{la:Gutman} implies that $a=3$. It follows that $G=C_{n,3}$, and so the 
theorem holds in Case 1. \\[1mm] 
{\sc Case 2:} $K$ is not $2$-connected. \\
We now bound each term on the right hand side of \eqref{eq:split-W-into-2} separately.  
Clearly, $K$ is an Eulerian graph of order $n-a+1$ but not a cycle.  
Since $G$ is a smallest counterexample to Theorem \ref{theo:Eulerian-with-2nd-largest-W}, the bound
in Theorem \ref{theo:Eulerian-with-2nd-largest-W} holds for $K$  unless $b <5$ or $b \in \{7,9\}$. 
However, since $b \geq \frac{n+1}{2}$ and $n\geq 26$, $b$ is not one of these exceptional
values and Theorem \ref{theo:Eulerian-with-2nd-largest-W} holds for $K$. Therefore, 
\begin{equation}   \label{eq:separable-bound-for-K} 
W(K) \leq W(C_{n-a+1,3}) \leq \frac{1}{8}(n-a+1)^3 - \frac{1}{4} (n-a+1)^2 + \frac{3}{2}(n-a+1) - 2. 
\end{equation}
It follows from Theorem \ref{theo:plesnik}(c) that
\begin{equation}   \label{eq:separable-bound-for-sigma(v)-in-K} 
\sigma_K(v) \leq \frac{1}{3}(n-a+1)(n-a). 
\end{equation}
As in Case 1, Theorem \ref{theo:plesnik} yields the following bounds for $W(H)$ and $\sigma_H(v)$  
\begin{equation}  \label{eq:separable-bound-for-H} 
W(H) \leq W(C_a) \leq \frac{1}{3}a^3 \ \ \textrm{and} \ \  
\sigma_H(v)  \leq \frac{1}{4}a^2. 
\end{equation} 
Substituting \eqref{eq:separable-bound-for-H}, \eqref{eq:separable-bound-for-K} and 
\eqref{eq:separable-bound-for-sigma(v)-in-K} 
into \eqref{eq:split-W-into-2} yields that 
\begin{eqnarray*}
W(G) & \leq & 
\frac{1}{8}a^3   
+ \frac{1}{8}(n-a+1)^3 - \frac{1}{4} (n-a+1)^2 + \frac{3}{2}(n-a+1) - 2 \\
& & +  \frac{1}{3}(a-1)(n-a+1)(n-a)
+ \frac{1}{4}(n-a)a^2.
\end{eqnarray*}
From equation \eqref{eq:value-W(C(n,3))} and the remark following it, we have 
\[ W(C_{n,3}) \geq \frac{1}{8}n^3 - \frac{1}{4}n^2 + \frac{11}{8}n - \frac{9}{4}. \]
Subtracting these two bounds we obtain, after simplification,
\[
W(C_{n,3}) - W(G)  \geq    \frac{1}{24} \Big\{ (a-1)n^2 + (a^2-18a+8)n
-2a^3 + 13 a^2  + 25 a - 39 \Big\}. 
\]
Denote the right hand side of the above inequality by $f(n,a)$. To complete the proof of the 
Lemma it suffices to show that $f(n,a)>0$ for $n\geq 26$ and $3 \leq a \leq \frac{n+1}{2}$. 
Now 
$\frac{\partial f}{\partial a}
 = \frac{1}{24} \big\{ n^2 + (2a-18)n  -6 a^2 + 26a  + 25  \big\}$. For constant $n$, this 
is a quadratic function of $a$ which is concave down and thus it attains its minimum for 
$a \in [3,\frac{n+1}{2}]$ at $a=3$ or $a= \frac{n+1}{2}$. Since 
for $a=3$ we have $\frac{\partial f}{\partial a} = \frac{1}{24}(n^2-12n+49) >0$, and 
for $a=\frac{n+1}{2}$ we have $\frac{\partial f}{\partial a} = \frac{1}{48}(n^2-14n+73) >0$, 
the derivative $\frac{\partial f}{\partial a}$ is positive for $3 \leq a \leq \frac{n+1}{2}$. 
It follows that the function $f$ is increasing in $a$ for constant $n$, and thus 
\[
W(C_{n,3}) - W(G)  \geq f(3)  = \frac{1}{24}\big( 2n^2 -37n +99\big),
\]
which is greater than $0$ for $n\geq 26$. This completes the proof of 
Theorem \ref{theo:Eulerian-with-2nd-largest-W}
\hfill$\Box$

\section{Eulerian Graphs with Small Wiener Index}

A natural question that arises in the context of of the Wiener index of 
Eulerian graphs is how small the Wiener index of an Eulerian graph can be.
For Eulerian graphs of given order, this was answered in \cite{GutCruRad2014}. 

\begin{prop} [Gutman, Cruz and Rada \cite{GutCruRad2014}]
Let $G$ be an Eulerian graph of order $n$, where $n\geq 3$. Then
\[ W(G) \geq \left\{ \begin{array}{cc}
       \binom{n}{2} & \textrm{if $n$ is odd,}\\
       \binom{n}{2} + \frac{n}{2} & \textrm{if $n$ is even.}
          \end{array}  \right.   \]
Equality holds if and only if $G$ is complete (for odd $n$), or $G$
is obtained from the complete graph by removing the edges of a perfect
matching (for even $n$). 
\end{prop}

Finding the minimum value of the Wiener index of Eulerian graphs 
becomes more challenging if not only the order, but also the size
of the graph is considered. We have the following lower bound on 
the Wiener index due to Plesn\'\i k \cite{Ple1984}.

\begin{prop}  \label{prop:lower-bound-on-W-for-given-size}
Let $G$ be a connected graph with $n$ vertices and $m$ edges. Then
\[ W(G) \geq 2 \binom{n}{2} -m. \]
Equality holds if and only if the diameter of $G$ is at most $2$. 
\end{prop}

Proposition \ref{prop:lower-bound-on-W-for-given-size} yields a lower bound on
the Wiener index of Eulerian graphs of given order and size. However, if $m$ is so small relative to $n$
that there is no Eulerian graph of diameter two of order $n$ and size $m$, 
then this bound is not sharp. 
The following result determines the minimum size of an Eulerian graph 
of order $n$ and diameter $2$. In the proof we use the fact that 
the minimum size of a graph of order $n$ and diameter $2$ not containing a vertex of 
degree $n-1$ is $2n-5$, which was proved by Erd\capitalhungarumlaut{o}s and R\'{e}nyi
\cite{ErdRen1962}, see also  \cite{HenSou2015}.

\begin{prop}
Let $G$ be an Eulerian graph of order $n$ and diameter two. Then
\[ m(G) \geq \left\{ \begin{array}{cc}
\frac{3}{2}(n-1) & \textrm{if $n$ is odd,} \\
2n-5 & \textrm{if $n$ is even.} 
\end{array} \right. \]
This bound is sharp for $n\geq 9$. 
\end{prop}

{\bf Proof:} First let $n$ be even. Since $G$ contains only vertices of even degree,
$G$ has no vertex of degree $n-1$. The above-mentioned result by 
Erd\capitalhungarumlaut{o}s and R\'{e}nyi
\cite{ErdRen1962} now proves that $m(G)\geq 2n-5$, as desired. \\
To see that the bound is sharp consider the graph obtained from a triangle with 
vertices $a$, $b$ and $c$ and a star $K_{1,n-4}$ by joining two of the leaves of 
the star to $a$, joining two other leaves to $b$, and joining the remaining 
$n-8$ leaves to $c$. (We note that this graph was already described in \cite{HenSou2015}.) \\
Now let $n$ be odd. If $G$ contains no vertex of degree $n-1$, then we have  $m\geq 2n-5$
as above, and the result follows. If $G$ contains a vertex of degree $n-1$, then all
other vertices have degree at least $2$, and so the degree sum of $G$ is at least 
$n-1 + 2(n-1)$, and so $m\geq \frac{3}{2}(n-1)$, as desired. \\
The graph obtained from  $\frac{n-1}{2}$ copies of the graph $K_2$ by adding an new vertex 
and joining it to each of the $n-1$ vertices shows that the bound is sharp. 
\hfill $\Box$ \\

This leads naturally to the following question which we pose as a problem.

\begin{que}
Given $n$ with $n\geq 9$, and $m$ with $m < 2n-5$ if $n$ is even and $m < \frac{3}{2}(n-1)$
if $n$ is odd. What is the minimum Wiener index of an Eulerian graph of order $n$ and
size $m$ and which graphs attain it?
\end{que}

\end{document}